\documentclass{kapproc}
\usepackage{latexsym,amssymb,graphicx,amsmath}

\upperandlowercase
\normallatexbib
\let\footnote\savefootnote

\begin{document}

\articletitle[A volume formula for generalized hyperbolic tetrahedra]{%
A volume formula for\\
generalized hyperbolic\\
tetrahedra\thanks{This research was partially supported by the Inamori Foundation.}}

\author{Akira {\sc Ushijima}}

\affil{Department of Mathematics,\\
Faculty of Science,\\
Kanazawa University,\\
920-1192 
JAPAN\footnote{Current address: Mathematics Institute, University of Warwick, Coventry CV4 7AL, UK}}
\email{ushijima@kenroku.kanazawa-u.ac.jp\footnote{Current e-mail address: 
{\emailfont ushijima@maths.warwick.ac.uk}}}

\makeatletter
        \renewcommand{\theequation}{%
                \thesection.\arabic{equation}}
        \@addtoreset{equation}{section}
\makeatother

\newtheorem{df}{Definition}[section]
\newtheorem{thm}[df]{Theorem}
\newtheorem{prop}[df]{Proposition}
\newtheorem{pbm}{Problem}

\newenvironment{abenumerate}{
 \begin{enumerate}
 \renewcommand{\labelenumi}{
\rm (\arabic{enumi})}
}{\end{enumerate}}

\newcommand{\bsquare}{\hbox{\rule{6pt}{6pt}}}
\newcommand{\Aut}{\mathop{\mathrm{Aut}}\nolimits}
\newcommand{\sgn}{\mathop{\mathrm{sgn}}\nolimits}
\newcommand{\ext}{\mathop{\mathrm{Ext}}\nolimits}
\newcommand{\vol}{\mathop{\mathrm{Vol}}\nolimits}
\newcommand{\li}{\mathop{\mathrm{Li_2}}\nolimits}

\begin{abstract}
A generalized hyperbolic tetrahedra is a polyhedron (possibly non-compact) \sloppy
with finite volume in hyperbolic space,
obtained from a tetrahedron by the polar truncation at the vertices lying outside the space.
In this paper it is proved that a volume formula for ordinary hyperbolic tetrahedra 
devised by J.~Murakami and M.~Yano can be applied to such ones.
There are two key tools for the proof;
one is so-called Schl\"afli's differential formula for hyperbolic polyhedra,
and the other is a necessary and sufficient condition for given numbers to be the dihedral angles of
a generalized hyperbolic simplex
with respect to their dihedral angles.
\end{abstract}

\begin{keywords}
hyperbolic tetrahedron,  Gram matrix, volume formula.
\end{keywords}

\begin{flushleft}
\begin{small}
{\bf 2000 Mathematics Subject Classifications:}  
Primary: 52A38; secondary: 51M09.
\end{small}
\end{flushleft}


\section{Introduction}
\label{sec_intro}

Obtaining volume formulae for basic polyhedra is
one of the basic and important problems in geometry.
In hyperbolic space this problem has been attacked from its birth.
An {\em orthoscheme\/} is a simplex of a generalization (with respect to dimensions) of a right triangle.
Since any hyperbolic polyhedra can be decomposed into finite number of orthoschemes, 
the first step of solving the problem is to obtain a volume formula for orthoschemes.
In three-dimensional case, 
N.~I.~Lobachevsky \cite{lo} found in 1839 a volume formula for orthoschemes.
Thus the next step is to find a formula for ordinary tetrahedra.

Although the lengths of six edges of tetrahedra determine their sizes and shapes,
in hyperbolic space it is known that the six dihedral angles also plays the same role.
Actually the valuables of the Lobachevsky's formula mentioned above are dihedral angles,
and following this way in this paper we obtain a volume formula
for {\em generalized\/} hyperbolic tetrahedra
in terms of the six dihedral angles
(the meaning of the term ``generalized" will be explained later).
Here we note that, by the fact mentioned previous paragraph,
we can construct such a formula from that of orthoschemes.
But in this case we have to calculate the dihedral angles of orthoschemes from those of tetrahedra.
It is emphasized that what we want to obtain is the formula
being able to calculate directly from the dihedral angles  of a given generalized tetrahedron.

The first answer to this problem was (at least the author is aware) given by W.-Y. Hsiang in 1988.
In \cite{hs} the formula was presented through an integral expression.
On the other hand, many volume formulae for hyperbolic polyhedra is presented
via Lobachevsky function or the dilogarithm function
(see, for example, \cite{ve}).
The first such presentation was given by Y.~Cho and H.~Kim in 1999 (see \cite{ck}).
Their formula was derived from that of ideal tetrahedra,
and has the following particular property:
due to the way of the proof it is not symmetric with respect to the dihedral angles.
Thus the next problem is to obtain its essentially symmetrical expression.
Such a formula was obtained by J.~Murakami and M.~Yano in about 2001.
In \cite{my} they derived it from the quantum $6 j$-symbol, but the proof is to reduce it to Cho-Kim's one.

\vspace{\bigskipamount}

In hyperbolic space it is possible to consider that
vertices of a polyhedron lie ``outside the space and its sphere at infinity."
For each such vertex there is a canonical way to cut off the vertex with its neighborhood.
This operation is called a {\em truncation\/} at the vertex, 
which will be defined more precisely in Section~\ref{sec_exist},
and we thus obtain convex polyhedron (possibly non-compact) with finite volume.
Such a polyhedron appears, for example, 
as fundamental polyhedra for hyperbolic Coxeter groups 
or building blocks of three-dimensional hyperbolic manifolds with totally geodesic boundary.
Since the volume of hyperbolic manifolds are topological invariants,
it is meaningful, also for $3$-manifold theory, to obtain a volume formula for such polyhedra.

R.~Kellerhals presented in 1989 that 
the formula for orthoschemes can be applied, without any modification,
to ``mildly" truncated ones at {\em principal vertices\/} 
(see Appendix
for the definition and \cite{ke} for the formula).
This result inspires that Murakami-Yano's formula also
may be applied to {\em generalized hyperbolic tetrahedra},
tetrahedra of each vertices being finite, ideal or truncated
(see Definition~\ref{df_trunc} for precise definition),
and this is what we will prove in this paper (see Theorem~\ref{thm_volumeformula}).
Here we note that, 
although the lengths of the edges emanating from a vertex at infinity are infinity,
the dihedral angles may be finite.
This is another reason why we take not the edge lengths but the dihedral angles
as the valuables of our volume formula.

The key tool for the proof is so-called Schl\"afli's differential formula,
a simple description for the volume differential of polyhedra 
as a function of the dihedral angles and the volume of the apices.
Since the volume formula is a function of the dihedral angles,
to apply Schl\"afli's differential formula 
we have to translate
the dihedral angles of a generalized tetrahedron to the lengths of its edges.
This requirement yields a necessary and sufficient condition of a set of positive numbers
to be the dihedral angles of a generalized simplex (see Theorem~\ref{thm_exist}).

\subsubsection*{Volume formula for generalized hyperbolic tetrahedra}
\label{subsubsec_volumeformula}

\begin{figure}[ht]
        \begin{center}
        \includegraphics[clip]{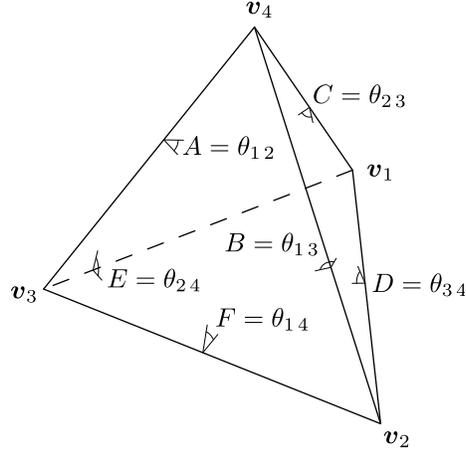}  
        \end{center}
        \caption{The dihedral angles of $T$}
        \label{fig_tetraangle}
\end{figure}

Let $T=T ( A,B,C,D,E,F )$ be a generalized tetrahedron in the three-dimensional
hyperbolic space ${\Bbb H}^{3}$ whose dihedral angles are $A,B,C,D,E,F$.
Here the configuration of the dihedral angles are as follows (see also Figure~\ref{fig_tetraangle}):
three edges corresponding to $A$, $B$ and $C$ arise from a vertex,
and the angle $D$ (resp.\ $E$, $F$) is put on the edge
opposite to that of $A$ (resp.\ $B$, $C$).
Let $G$ be the {\em Gram matrix\/} of $T$ defined as follows:
$$
G := \left(
\begin{array}{@{\,}cccc@{\,}}
1 & -\cos A & -\cos B & -\cos F \\
-\cos A & 1 & -\cos C & -\cos E \\
-\cos B & -\cos C & 1 & -\cos D \\
-\cos F & -\cos E & -\cos D & 1
\end{array} 
\right) .
$$

Let $a := \exp \sqrt{-1} A, b := \exp \sqrt{-1} B, \ldots, f := \exp \sqrt{-1} F$, 
and let $U ( z, T )$ be the complex valued function defined as follows:
\begin{eqnarray*}
U ( z, T ) &:=& \frac{1}{2} \left\{ \li (z) + \li (a b d e z) + \li (a c d f z) + \li (b c e f z) \right.\\
&& \text{ } \left. - \li (- a b c z) - \li (- a e f z) - \li (- b d f z) - \li (- c d e z) \right\},
\end{eqnarray*}
where $\li (z)$ is the dilogarithm function defined by the analytic continuation of the following integral:
$$
\li (x) := - \int_{0}^{x} \frac{\log (1-t )}{t} \mathrm{d} t \quad \text{ for a positive real number $x$.}
$$
We denote by $z_1$ and $z_2$ the two complex numbers defined as follows:
\begin{eqnarray*}
z_1 &:=& -2 \frac{\sin A \sin D + \sin B \sin E +\sin C \sin F - \sqrt{\det G}}
{a d + b e + c f + a b f + a c e + b c d + d e f + a b c d e f} , \\
z_2 &:=& -2 \frac{\sin A \sin D + \sin B \sin E +\sin C \sin F + \sqrt{\det G}}
{a d + b e + c f + a b f + a c e + b c d + d e f + a b c d e f} .
\end{eqnarray*}
%
%
%
%
%
%
%
%
\begin{thm}
\label{thm_volumeformula}
The volume $\vol (T)$ of a generalized tetrahedron $T$ is given as follows:
\begin{equation}
\vol (T) = \frac{1}{2} \Im \left( U ( z_1, T ) - U ( z_2, T ) \right) , 
\label{eqn_volumeformula}
\end{equation}
where $\Im$ means the imaginary part.
\end{thm}

\begin{acknowledgments}
The author would like to give his thanks to Professor Jun Murakami and
Mr Masakazu Yano for useful discussion about their results. 
The author also would like to give his sincere gratitude to the referee for his careful reading and advices.
\end{acknowledgments}

\section{Preliminaries}
\label{sec_pre}

In this section we review several well-known facts about hyperbolic geometry.
See, for example, \cite{us} for more precise explanation and the proofs of the propositions.

The {\em $n+1$-dimensional Lorentzian space\/}  
${\Bbb E}^{1,n}$
is the real vector space 
${\Bbb R}^{n+1}$ of  
dimension $n+1$ with the 
{\em Lorentzian inner product\/} 
$ \left \langle \text{\boldmath $x$}, 
\text{\boldmath $y$} \right \rangle 
:= - x_0 y_0 + x_1  y_1 + \cdots + x_n y_n$, 
where $\text{\boldmath $x$}
= \left( x_0, x_1, \ldots , x_n \right)$ and 
$\text{\boldmath $y$}
= \left( y_0, y_1, \ldots , y_n \right)$. 
Let \sloppy
$H_T := \left\{ \left.  \text{\boldmath $x$} \in {\Bbb E}^{1,n} \, \right|
\left \langle \text{\boldmath $x$}, \text{\boldmath $x$} \right \rangle = - \, 1 \, \right\}$ 
be the (standard) hyperboloid of two sheets, and let
$H_T^{+} := \left\{ \left. 
\text{\boldmath $x$} \in
{\Bbb E}^{1,n} \, \right|
\left \langle \text{\boldmath $x$}, 
\text{\boldmath $x$} \right \rangle 
= - \, 1 \text{ and } x_0 > 0 \, \right\}$ 
be its upper sheet.
The restriction of the quadratic form 
induced by 
$\left \langle \, \cdot ,  \cdot \, \right \rangle$ 
on ${\Bbb E}^{1,n}$ to the tangent spaces of 
$H_T^{+}$ is positive definite and gives 
a Riemannian metric on $H_T^{+}$. 
The space obtained from $H_T^{+}$ equipped with 
the metric above is called the
{\em hyperboloid model\/}  
of the $n$-dimensional hyperbolic space, and we denote it by ${\Bbb H}^n$. 
Under this metric the hyperbolic distance, say $d$, 
between two points {\boldmath $x$} and {\boldmath $y$} can be measured by the following formula:
%
%
%
%
\begin{equation}
\label{eqn_dist}
\left \langle \text{\boldmath $x$}, 
\text{\boldmath $y$} \right \rangle 
= - \cosh d \, .
\end{equation}
Let $L := \left\{ \left. \text{\boldmath $x$} \in {\Bbb E}^{1,n} \, \right|
\left \langle \text{\boldmath $x$}, \text{\boldmath $x$} \right \rangle = 0 \, \right\}$
be the (standard) cone and let
$L^{+} := \left\{ \left. 
\text{\boldmath $x$} \in
{\Bbb E}^{1,n} \, \right|
\left \langle \text{\boldmath $x$}, 
\text{\boldmath $x$} \right \rangle 
= 0 \text{ and } x_0 > 0 \, \right\}$ be its upper half. \sloppy
Then
a ray in $L^{+}$ 
started from the origin {\boldmath $o$} 
corresponds to a point in the
ideal boundary of ${\Bbb H}^n$. 
The set of such rays forms 
the sphere at infinity,  
and we denote it by ${\bf S}^{n-1}_\infty$.
Then each ray in $L^{+}$ becomes a point at
infinity of ${\Bbb H}^n$. 

Let us denote by $\cal P$ 
the radial projection from 
${\Bbb E}^{1,n} - 
\left\{ \left. \text{\boldmath $x$} 
\in {\Bbb E}^{1,n} \, \right|
x_0 = 0 \, \right\}$ 
to an affine hyperplane 
${\Bbb P}_1^{n} := 
\left\{ \left. \text{\boldmath $x$} 
\in {\Bbb E}^{1,n} \, \right|
x_0 = 1 \, \right\}$
along the ray from the origin {\boldmath $o$}. 
The projection $\cal P$ 
is a homeomorphism 
on ${\Bbb H}^n$ to the $n$-dimensional 
open unit ball ${\bf B}^n$ in ${\Bbb P}_1^n$ 
centered at 
the origin 
$\left( 1, 0, 0, \ldots , 0 \right)$  
of ${\Bbb P}_1^n$, 
which gives the {\em projective model\/} 
of ${\Bbb H}^n$. 
The affine hyperplane ${\Bbb P}_1^n$ contains 
not only ${\bf B}^n$ and its set theoretic 
boundary $\partial {\bf B}^n$ 
in ${\Bbb P}_1^n$, 
which is canonically identified with 
${\bf S}^{n-1}_\infty$, 
but also the outside of 
the compactified projective model 
$\overline{{\bf B}^n} := 
{\bf B}^n \sqcup \partial {\bf B}^n 
\approx {\Bbb H}^n \sqcup {\bf S}^{n-1}_\infty$. 
So $\cal P$ can be naturally extended to the mapping
from ${\Bbb E}^{1,n} - 
\left\{ \text{\boldmath $o$} \right\}$  
to the 
$n$-dimensional real projective space 
${\Bbb P}^n := {\Bbb P}_1^n 
\sqcup {\Bbb P}_{\infty}^n$, 
where ${\Bbb P}_{\infty}^n$
is the set of lines in the 
affine hyperplane 
$\left\{ \left. \text{\boldmath $x$} \in 
{\Bbb E}^{1,n} \, \right| x_0 = 0 \, \right\}$ 
through {\boldmath $o$}. 
We denote by 
$\ext \overline{{\bf B}^n}$
the exterior of 
$\overline{{\bf B}^n}$ in ${\Bbb P}^n$. 

The (standard) {\em hyperboloid of one sheet\/} 
$H_S$
is defined to be \sloppy
$H_S := \left\{ \left. \text{\boldmath $x$} 
\in {\Bbb E}^{1,n} \, \right|
\left \langle \text{\boldmath $x$}, 
\text{\boldmath $x$} \right \rangle  
= 1 \, \right\}$. 
For an arbitrary point {\boldmath $u$} 
in $H_S$, 
we define a half-space 
$R_{\text{\boldmath $u$}}$ and a hyperplane 
$P_{\text{\boldmath $u$}}$ in ${\Bbb E}^{1,n}$ 
as follows: 
\begin{eqnarray*}
R_{\text{\boldmath $u$}} 
&:=& 
\left\{ \text{\boldmath $x$} \in
{\Bbb E}^{1,n} \left| \, 
\left \langle \text{\boldmath $x$}, 
\text{\boldmath $u$}  \right \rangle \leq 0
\right. \right\} \, , \\
&&\\
P_{\text{\boldmath $u$}} 
&:=& 
\left\{ \text{\boldmath $x$} \in
{\Bbb E}^{1,n} \left| \, 
\left \langle \text{\boldmath $x$}, 
\text{\boldmath $u$} \right \rangle = 0
\right. \right\} = 
\partial R_{\text{\boldmath $u$}} \, .
\end{eqnarray*}
We denote by $\Gamma_{\text{\boldmath $u$}}$ 
(resp.\ $\Pi_{\text{\boldmath $u$}}$)
the intersection of $R_{\text{\boldmath $u$}}$ 
(resp.\ $P_{\text{\boldmath $u$}}$) and 
${\bf B}^n$.
Then $\Pi_{\text{\boldmath $u$}}$ is a geodesic hyperplane in ${\Bbb H}^{n}$,
and the correspondence between the points in $H_S$
and the half-spaces $\Gamma_{\text{\boldmath $u$}}$ in ${\Bbb H}^{n}$
is bijective.
We call {\boldmath $u$} 
a {\em normal vector\/} to 
$P_{\text{\boldmath $u$}}$ 
(or $\Pi_{\text{\boldmath $u$}}$). 
The following two propositions are on relationships
between the Lorentzian inner product and geometric objects.
%
%
%
%
%
%
%
%
\begin{prop}\label{prop_linear}
Let\/ {\boldmath $x$} and\/ {\boldmath $y$} 
be arbitrary two non-parallel points in\/ $H_S$.
Then one of the followings hold:
\begin{abenumerate}
\item \label{prop_linear_intersect}
Two geodesic hyperplanes\/ $\Pi_{\text{\boldmath $x$}}$ and\/ $\Pi_{\text{\boldmath $y$}}$
intersect
if and only if\/ 
$\left| \left \langle \text{\boldmath $x$}, \text{\boldmath $y$} \right \rangle \right| < 1$.
In this case the (hyperbolic) angle, say $\theta$,
between them
measured in\/ $\Gamma_{\text{\boldmath $x$}}$ and\/ $\Gamma_{\text{\boldmath $y$}}$
is calculated by the following formula:
\begin{equation}
\label{eqn_angle}
\left \langle \text{\boldmath $x$}, \text{\boldmath $y$} \right \rangle = - \cos \theta \, .
\end{equation}
\item \label{prop_linear_ultraparallel}
Two geodesic hyperplanes\/ $\Pi_{\text{\boldmath $x$}}$ and\/ $\Pi_{\text{\boldmath $y$}}$
never intersect in\/ $\overline{{\bf B}^n}$, i.e., they intersect in\/ $\ext \overline{{\bf B}^n}$,
if and only if\/ 
$\left| \left \langle \text{\boldmath $x$}, \text{\boldmath $y$} \right \rangle \right| > 1$.
In this case the (hyperbolic) distance, say $d$,
between them
is calculated by
\begin{equation}
\label{eqn_distance}
\left| \left \langle \text{\boldmath $x$}, \text{\boldmath $y$} \right \rangle \right| = \cosh d \, ,
\end{equation}
and then\/ $\Pi_{\text{\boldmath $x$}}$ and\/ $\Pi_{\text{\boldmath $y$}}$
are said to be\/ {\em ultraparallel}.
\item \label{prop_linear_parallel}
Two geodesic hyperplanes\/ $\Pi_{\text{\boldmath $x$}}$ and\/ $\Pi_{\text{\boldmath $y$}}$
intersect not in\/ ${\bf B}^n$ but in\/ $\partial {\bf B}^n$
if and only if\/ 
$\left| \left \langle \text{\boldmath $x$}, \text{\boldmath $y$} \right \rangle \right| = 1$.
In this case the angle and the distance between them is $0$,
and then\/ $\Pi_{\text{\boldmath $x$}}$ and\/ $\Pi_{\text{\boldmath $y$}}$
are said to be\/ {\em parallel}.
\end{abenumerate}
\end{prop}
%
%
%
%
%
%
%
%
\begin{prop}\label{prop_distPoPl}
Let\/ {\boldmath $x$} be a point in\/ ${\bf B}^n$ 
and let\/ $\Pi_{\text{\boldmath $y$}}$ be a geodesic hyperplane 
whose normal vector is\/ $\text{\boldmath $y$} \in H_S$ with 
$\left \langle \text{\boldmath $x$}, \text{\boldmath $y$} \right \rangle < 0$.
Then the distance $d$ between\/  {\boldmath $x$} and\/ $\Pi_{\text{\boldmath $y$}}$ 
is obtained by the following formula:
\begin{equation}\label{eqn_distPoPl}
\left \langle \text{\boldmath $x$}, \text{\boldmath $y$} \right \rangle =  - \sinh d \, .
\end{equation}
\end{prop}

Let {\boldmath $v$} be a point in $\ext \overline{{\bf B}^n}$.
Then ${\cal P}^{-1} (\text{\boldmath $v$}) \cap H_S$ consists of two points
and, independent of the choice of
$\widetilde{\text{\boldmath $v$}} \in {\cal P}^{-1} (\text{\boldmath $v$}) \cap H_S$,
we can define the same hyperplane $\Pi_{\widetilde{\text{\boldmath $v$}}}$.
We call $\Pi_{\widetilde{\text{\boldmath $v$}}}$ 
the {\em polar geodesic hyperplane to} {\boldmath $v$}
and {\boldmath $v$} the {\em pole\/} of $\Pi_{\widetilde{\text{\boldmath $v$}}}$.
Then the following proposition holds:
%
%
%
%
%
%
%
%
\begin{prop}
\label{prop_tangent}
Let\/ {\boldmath $v$} be a point in $\ext \overline{{\bf B}^n}$.
\begin{abenumerate}
\item Any hyperplane through {\boldmath $v$} with intersecting ${\bf B}^n$
is perpendicular to $\Pi_{\widetilde{\text{\boldmath $v$}}}$ in ${\Bbb H}^n$.
\item Let\/ {\boldmath $u$} be a limit point of $\Pi_{\widetilde{\text{\boldmath $v$}}}$,
i.e., $\text{\boldmath $u$} \in P_{\text{\boldmath $v$}} \cap \partial {\bf B}^n$.
Then the line through {\boldmath $u$} and {\boldmath $v$} is tangent to $\partial {\bf B}^n$.
\end{abenumerate}
\end{prop}

\section{A characterization theorem for the dihedral angles of generalized hyperbolic simplices}
\label{sec_exist}

\begin{figure}[ht]
        \begin{center}
        \includegraphics[clip]{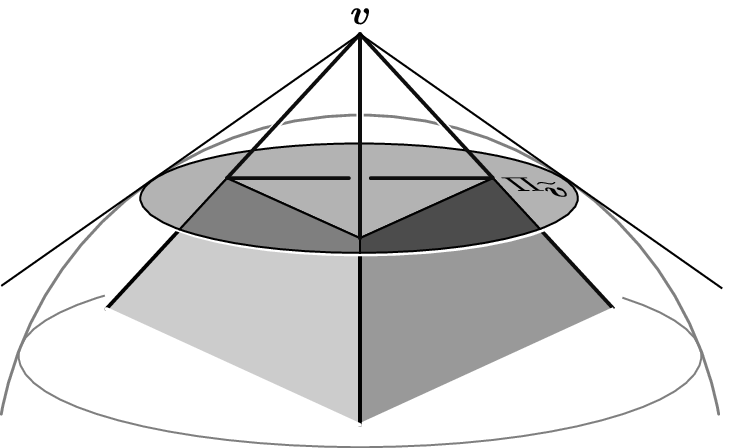}
        \end{center}
        \caption{The truncation of $\triangle$ at {\boldmath $v$}}
        \label{fig_truncation}
\end{figure}

Unlike spherical or Euclidean one, 
in hyperbolic geometry 
we can consider not only points at infinity but also points ``beyond infinity."
This situation can be easily seen in the projective ball model,
and it extends 
the concept of simplices.
Such simplices are called the \em{generalized simplices}\rm{}.
We start this section with their precise definition.
In this section we assume $n \geq 3$.
%
%
%
%
%
%
%
%
\begin{df} 
\label{df_trunc}
\em{
Let $\triangle$ be a simplex in ${\Bbb P}_{1}^{n}$ containing the origin of ${\Bbb P}_{1}^{n}$ 
(remark; 
since any polyhedron in $\overline{{\bf B}^n}$
can be moved to the one containing the origin
by the action of orientation-preserving isomorphisms of ${\Bbb H}^{n}$,
this assumption does not lose generality in hyperbolic geometry).
Suppose each of the interior of its ridge (i.e., $(n - 2)$-dimensional face) 
intersects $\overline{{\bf B}^n}$.
\begin{enumerate}
\item 
Let {\boldmath $v$} be a vertex of $\triangle$ in $\ext \overline{{\bf B}^n}$.
The \em{(polar) truncation\/}\rm{} of $\triangle$ at {\boldmath $v$} is an operation of 
omitting the pyramid of apex {\boldmath $v$} 
with base $(n - 1)$-dimensional simplex
$\Pi_{\widetilde{\text{\boldmath $v$}}} \cap \triangle$,
and then capping the open end by 
$\Pi_{\widetilde{\text{\boldmath $v$}}} \cap \triangle$ (see Figure~\ref{fig_truncation}).
\item 
The \em{truncated simplex\/}\rm{}, say $\triangle'$, 
is a polyhedron in $\overline{{\bf B}^n}$ obtained from $\triangle$ 
by the truncation at all vertices in $\ext \overline{{\bf B}^n}$.
\item
A \em{generalized simplex\/}\rm{} in ${\Bbb H}^{n}$ is 
a polyhedron which is either a simplex in the ordinary sense or a truncated simplex 
described above.
\end{enumerate}
}\end{df} \vspace{\medskipamount}

In this paper we will regard vertices of $\triangle$ as those of $\triangle'$.
This means, for example,
a vertex of $\triangle$ in $\ext \overline{{\bf B}^n}$ is no longer that of $\triangle'$
in the ordinary sense, since it is not an element of $\triangle'$,
but we will call it a ``vertex of $\triangle'$."
Conversely vertices (in the ordinary sense) arisen by the truncation will not be mentioned 
``vertices of $\triangle'$" in this paper.
Under this usage,
we call a vertex {\boldmath $v$} of $\triangle'$ \em{finite\/}\rm{} 
(resp.\ \em{ideal}\rm{}, \em{ultraideal\/}\rm{})
when $\text{\boldmath $v$} \in {\bf B}^n$ 
(resp.\ $\partial {\bf B}^n, \ext \overline{{\bf B}^n}$).

Let $\sigma^n$ be an $n$-dimensional generalized simplex in the projective model
with vertices $\left\{ \text{\boldmath $v$}_i \right\}_{i \in I}$,
where $I := \left\{ 1,2, \ldots , n+1 \right\}$. 
Then we regard the lift of its vertices to ${\Bbb E}^{1,n}$ as follows:
If a vertex {\boldmath $v$} is finite, 
then the lift is uniquely determined by ${\cal P}^{-1} (\text{\boldmath $v$}) \cap H_T^+$.
If a vertex is ultraideal, there are two choices of the lift,
and we choose the one defining the half-space containing $\sigma^n$.
If a vertex is ideal, we do not need at this point to determine the exact lift in $L^{+}$.
The \em{$i$-th facet\/}\rm{} of $\sigma^n$ is the ($n-1$)-dimensional face of 
$\sigma^n$ opposite to $\text{\boldmath $v$}_i$.

\vspace{\bigskipamount}

Let $G$ be a matrix of order $n+1$.
We prepare several notations on matrices which will be used later.
\begin{enumerate}
 \item We denote by $G_{i \, j}$ the submatrix of order $n$ obtained from $G$ 
 by deleting its $i$th row and $j$th column.
 \item We denote a cofactor of $G$ by $c_{i \, j} := (-1)^{i+j} \det G_{i \, j}$.
 \item Suppose $G$ is a real symmetric matrix. 
 Then we denote by $\sgn G$ the signature of $G$, i.e.,
 if $\sgn G = (a,b)$ then $G$ has $a$ positive and $b$ negative eigenvalues.
\end{enumerate}
The following theorem tells us a necessary and sufficient condition of a set of positive numbers
to be the dihedral angles of a generalized simplex.
Calculations given in the proof will be used in that of the main theorem.
%
%
%
%
%
%
%
%
\begin{thm}
\label{thm_exist}
Suppose the following set of positive numbers is given:
$$
\left\{ \left. \theta_{i \, j} \in \left[ 0, \pi \right] \, \right| 
i, j \in I,
\theta_{i \, j} = \theta_{j \, i},
\theta_{i \, j} = \pi \text{ {\rm iff} } i=j
\right\}.
$$
Then the following two conditions are equivalent:
\begin{abenumerate}
\item \label{exist}
There exists a generalized hyperbolic simplex in ${\Bbb H}^{n}$ 
and the dihedral angle between $i$-th facet and $j$-th facet is $\theta_{i \, j}$. 
\item \label{gram}
The real symmetric matrix $G := \left( -\cos \theta_{i \, j} \right)$ of order $n+1$
satisfies the following two conditions:
\begin{itemize}
\item[\rm (a)] $\sgn G = (n,1)$,
\item[\rm (b)] $c_{i \, j} > 0$ for any $i,j \in I$ with $i \not= j$.
\end{itemize}
\end{abenumerate}
\end{thm} \vspace{\bigskipamount}

\em{Proof.}\rm{}
First we assume (\ref{gram}) holds.
By (a) there exists an orthogonal matrix, say $U$, of order $n+1$ such that
the conjugate of $G$ by $U$ is as follows:
$$
G = {}^t \! U 
\left(
\begin{array}{@{\,}ccccc@{\,}}
-1 & 0 & 0 & \cdots & 0 \\
0 & 1 & 0 & \cdots & 0 \\
0 & 0 & 1 & \cdots & 0 \\
\vdots & \vdots & \vdots & \ddots & \vdots \\
0 & 0 & 0 & \cdots & 1
\end{array}
\right) U .
$$
Suppose $U = \left( \text{\boldmath $u$}_1, \text{\boldmath $u$}_2, \ldots , \text{\boldmath $u$}_{n+1} \right)$.
Since (a) holds, $\det U \not= 0$. So $\left\{ \text{\boldmath $u$}_1, \text{\boldmath $u$}_2, \ldots , \text{\boldmath $u$}_{n+1} \right\}$
are linearly independent and thus becomes a basis of ${\Bbb E}^{n,1}$.
Now we define a vector $\text{\boldmath $w$}_i$ as follows:
$$
\text{\boldmath $w$}_i := \sum_{k=1}^{n+1} c_{i \, k} \, \text{\boldmath $u$}_k . 
$$
Then it is easy to see that, for any $i, j \in I$, the equation
$\left< \text{\boldmath $w$}_i , \text{\boldmath $u$}_j \right> 
= \delta_{i \, j} \det G$
holds. Here $\delta_{i \, j}$ means the Kronecker delta.
Using this relation, we obtain
$\left< \text{\boldmath $w$}_i , \text{\boldmath $w$}_j \right> 
= c_{i \, j} \det G$.
Now we define $\text{\boldmath $v$}_i$ as follows:
\begin{equation}
\text{\boldmath $v$}_i := 
\left\{ \begin{array}{ll}
\displaystyle{\frac{\text{\boldmath $w$}_i}{\sqrt{\left| c_{i \, i} \det G \right|}}} & 
\text{if } c_{i \, i} \not= 0, \\ & \\
\text{\boldmath $w$}_i & \text{if } c_{i \, i} = 0.
\end{array} \right.
\label{eq_VertexPosition}
\end{equation}
The goal now is to show $\left\{ \text{\boldmath $v$}_1, \text{\boldmath $v$}_2, \ldots , \text{\boldmath $v$}_{n+1} \right\}$ constructs the generalized simplex in question.

Here we calculate $\left< \text{\boldmath $v$}_i , \text{\boldmath $v$}_j \right>$;
\begin{equation}
\left< \text{\boldmath $v$}_i , \text{\boldmath $v$}_j \right> = 
\left\{ \begin{array}{ll}
\displaystyle{\frac{-c_{i \, j}}{\sqrt{\left| c_{i \, i} c_{j \, j} \right|}}} & 
\text{if } c_{i \, i} c_{j \, j} \not= 0 , \\ & \\
-c_{i \, j} \displaystyle{\sqrt{\left|  \frac{\det G}{c_{i \, i}} \right|}} & 
\text{if } c_{i \, i} \not= 0 \text{ and } c_{j \, j} = 0 , \\ & \\
c_{i \, j} \det G & \text{if } c_{i \, i} = c_{j \, j} = 0 .
\end{array} \right. \label{eqn_vertices}
\end{equation}
When $c_{i \, i} > 0$,  
$\left< \text{\boldmath $v$}_i , \text{\boldmath $v$}_i \right> = -1$, i.e., 
$\text{\boldmath $v$}_i \in H_T$,
and when $c_{i \, i} < 0$,
$\left< \text{\boldmath $v$}_i , \text{\boldmath $v$}_i \right> = 1$, i.e., 
$\text{\boldmath $v$}_i \in H_S$.
Suppose $c_{i \, i} = 0$. 
Then, by (\ref{eqn_vertices}), $\left< \text{\boldmath $v$}_i , \text{\boldmath $v$}_i \right> = 0$.
Furthermore we have $\text{\boldmath $v$}_i = \text{\boldmath $w$}_i \not= \text{\boldmath $o$}$
since $\left\{ \text{\boldmath $u$}_1, \text{\boldmath $u$}_2, \ldots , \text{\boldmath $u$}_{n+1} \right\}$
are linearly independent and $c_{i \, j} \not= 0$ as $i \not= j$.
So $\text{\boldmath $v$}_i \in L - \left\{ \text{\boldmath $o$} \right\}$.

Previous calculations also show that
$\left< \text{\boldmath $v$}_i , \text{\boldmath $v$}_j \right> > 0$ 
when $i \not= j$.
So all $\text{\boldmath $v$}_i$'s 
contained in $H_T \sqcup \left( L - \left\{ \text{\boldmath $o$} \right\} \right)$
are in the upper or lower half space in ${\Bbb R}^{n+1}$.
Thus, exchanging all $\text{\boldmath $v$}_i$'s 
with $-\text{\boldmath $v$}_i$'s if necessary,
we can assume that
$\text{\boldmath $v$}_i \in H_T^{+} \sqcup L^{+} \sqcup H_S$.

By the construction of $\text{\boldmath $w$}_i$,
$\left\{ \text{\boldmath $w$}_1, \text{\boldmath $w$}_2, \ldots , \text{\boldmath $w$}_{n+1} \right\}$
are linearly independent, 
and so are $\left\{ \text{\boldmath $v$}_1, \text{\boldmath $v$}_2, \ldots , \text{\boldmath $v$}_{n+1} \right\}$.
Thus the convex hull $C$ of $\left\{ \text{\boldmath $v$}_1, \text{\boldmath $v$}_2, \ldots , \text{\boldmath $v$}_{n+1} \right\}$ is an $n$-simplex in ${\Bbb E}^{n+1}$.
Let $\sigma$ be the intersection of $C$ and 
$\bigcap_{\text{\boldmath $v$}_i \in H_S} R_{\text{\boldmath $v$}_i}$.
Then $P \left( \sigma \right)$ is the simplex in question and we obtain (\ref{exist}).

\vspace{\medskipamount}

Secondly we assume (\ref{exist}).
Let $\sigma$ be a generalized $n$-simplex with vertices 
$\left\{ \text{\boldmath $v$}_1, \text{\boldmath $v$}_2, \ldots , \text{\boldmath $v$}_{n+1} \right\}$.
We denote by $\text{\boldmath $u$}_i$ the outward unit normal vector 
to the facet opposite to $\text{\boldmath $v$}_i$.
Then, by the definition of the generalized simplex, 
$\text{\boldmath $u$}_i \in H_S$.
Since $- \cos \theta_{i \, j} = 
\left< \text{\boldmath $u$}_i , \text{\boldmath $u$}_j \right>$,
we have $\sgn G = \left( n,1 \right)$, which is the condition (a).

As before, we define a vector $\text{\boldmath $w$}_i$ as follows:
$$
\text{\boldmath $w$}_i := \sum_{k=1}^{n+1} c_{i \, k} \, \text{\boldmath $u$}_k . 
$$
By the definition of $\text{\boldmath $u$}_j$, 
$\left< \text{\boldmath $v$}_i , \text{\boldmath $u$}_j \right> = 0$
when $i \not= j$.
So $\text{\boldmath $v$}_i \in 
\bigcap_{i \not= j \in I} \text{\boldmath $u$}_{j}^{\bot}$.
On the other hand, $\text{\boldmath $w$}_i \in 
\bigcap_{i \not= j \in I} \text{\boldmath $u$}_{j}^{\bot}$
since $\left< \text{\boldmath $w$}_i , \text{\boldmath $u$}_j \right>
= \delta_{i \, j} \det G$.
The dimension of $\bigcap_{i \not= j \in I} \text{\boldmath $u$}_{j}^{\bot}$
is one, which implies that 
$\text{\boldmath $w$}_i$ is parallel to $\text{\boldmath $v$}_i$.
Especially $\text{\boldmath $w$}_i = d_i \text{\boldmath $v$}_i$
for some $d_i > 0$.
Thus we have
$$
\left< \text{\boldmath $w$}_i , \text{\boldmath $w$}_j \right>
= d_i d_j \left< \text{\boldmath $v$}_i , \text{\boldmath $v$}_j \right> 
= c_{i \, j} \det G . 
$$
So $c_{i \, j} = d_i d_j \left< \text{\boldmath $v$}_i , \text{\boldmath $v$}_j \right> / \det G$.
By the definition of generalized simplices together with the lifts of the vertices,
$\left< \text{\boldmath $v$}_i , \text{\boldmath $v$}_j \right> < 0$ for any $i \not= j \in I$.
Thus we have the condition (b), and the proof completes.
\hfill $\square$
\vspace{\bigskipamount}\\
{\bf Remarks.}
\begin{enumerate}
\item The proof of Theorem~\ref{thm_exist} follows that of THEOREM in \cite{lu},
where it is proved a necessary and sufficient condition 
for a given set of positive real numbers to be the dihedral angles 
of a hyperbolic simplex in the ordinary sense.
Moreover such conditions for spherical and Euclidean simplices are also presented in \cite{lu}.
\item We call the matrix $G$ appeared in the previous theorem the {\em Gram matrix}\/
of an $n$-dimensional generalized simplex $\sigma$.
Here it should be noted that this definition is slight different 
from the ordinary one (see, for example, \cite{vi}),
since we do not deal with the normal vectors of the faces obtained by truncation.
\item By the definition (\ref{eq_VertexPosition}), 
we can characterize
vertices from the Gram matrix;
if $c_{i \, i} > 0$ (resp.\ $c_{i \, i} = 0, c_{i \, i} < 0$), 
then the vertex opposite to the $i$-th facet is finite (resp.\ ideal, ultraideal).
\end{enumerate}

\section{Schl\"afli's differential formula and the maximal volume of generalized hyperbolic tetrahedra}
\label{sec_schlafli}

Now we go on to the
main purpose of this paper, 
namely obtaining a volume formula for generalized hyperbolic tetrahedra.
So, from now on, we only consider the three-dimensional case.
The key tool for a proof of our volume formula is so-called ``Schl\"afli's differential formula," 
which gives us a relationship between volume formulas and their dihedral angles.
In this section we present the formula and also give another usage of the differential formula, 
which is a determination of the maximal volume amongst the generalized hyperbolic tetrahedra.

In the three-dimensional case, Schl\"afli's differential formula is written as follows (see \cite[{\S}2]{ke}).
%
%
%
%
%
%
%
%
\begin{thm}
\label{thm_schlafli}
Let $T$ be a generalized tetrahedra in ${\Bbb H}^{3}$ with dihedral angles $\theta_i$
and corresponding edges of length $l_i$ for $i=1,2,\ldots,6$.
We denote by $\vol$ the volume function on the set of generalized tetrahedra.
Then the following equation holds:
$$
\mathrm{d} \vol (T) = - \frac{1}{2} \sum_{i=1}^{6} l_i \, \mathrm{d} \theta_i .
$$
\end{thm}

U.~Haagerup and H.~Munkholm proved in \cite{hm} that, 
in hyperbolic space of dimension greater than one, 
a simplex is of maximal volume if and only if it is ideal and regular\footnote{More generally, 
it is proved by N.~Peyerimhoff that the regular simplex is of maximal volume 
amongst all simplices in a closed geodesic ball in ${\Bbb H}^n$; see \cite[{\sc Theorem 1}]{pe}.}.
Once we extend the class of tetrahedra to the generalized ones in this paper,
it can be determined that, by using this formula,
the maximal volume generalized tetrahedron.
%
%
%
%
%
%
%
%
\begin{thm}
\label{thm_max}
The maximal volume amongst the generalized hyperbolic tetrahedra 
is realized by the ideal right-angled octahedron,
and the volume is $8 \Lambda(\pi/4)$,
where $\Lambda (x)$ is the function, called the\/ {\em Lobachevsky function}, defined as follows:
$$
\Lambda (x) := - \int_{0}^{x} \log \left| 2 \sin t \right| \mathrm{d} t .
$$
\end{thm}
\subsubsection*{Proof of Theorem~\ref{thm_max}.}
Since any tetrahedron is contained in an ideal hyperbolic tetrahedron, 
the maximal volume is attained by generalized tetrahedra whose vertices are outside the hyperbolic space.
Furthermore  Schl\"afli's differential formula says that 
the smaller the length of an internal edges is, the larger the volume of a generalized
tetrahedron is.
So the maximal volume is attained by a truncated tetrahedron with all internal edges having length $0$, 
and by the definition of the truncation, its dihedral angles are all equal to $\pi/2$.
Namely the right-angled ideal octahedron has a unique generalized tetrahedron with maximal volume.
The actual value of its volume will be calculated in the proof of the main theorem. \hfill $\square$

\section{Proof of a volume formula for generalized hyperbolic tetrahedra}
\label{sec_formula}

Let ${\cal V} (T) := \Im \left( U ( z_1, T ) - U ( z_2, T ) \right) / 2$ 
be the right hand side of the formula (\ref{eqn_volumeformula}).
Then,
by Theorem~\ref{thm_schlafli},
it is necessary to check that
the partial derivative of ${\cal V} (T)$
with respect to each dihedral angle
is equal to $-1/2$ times the length of the edge corresponding to it.
Since the formula is symmetric under permutation of the dihedral angles,
from now on we focus on a dihedral angle, say $A$.

\vspace{\bigskipamount}

We first calculate the length of the edge corresponding to $A$.
Let $\text{\boldmath $v$}_4$ be the vertex of being common endpoint 
of three edges corresponding to $A$, $B$ and $C$,
and let $\text{\boldmath $v$}_1$
(resp.\ $\text{\boldmath $v$}_2$, $\text{\boldmath $v$}_3$)
be the other endpoint of the edge corresponding to $C$ (resp.\ $B$, $A$).
It should be noted that we also regard $\text{\boldmath $v$}_i$ as the vector in ${\Bbb E}^{1,3}$
corresponding the vertex.
Recall that $G$ is the Gram matrix of $T$ defined as follows:
$$
G := \left(
\begin{array}{@{\,}cccc@{\,}}
1 & -\cos A & -\cos B & -\cos F \\
-\cos A & 1 & -\cos C & -\cos E \\
-\cos B & -\cos C & 1 & -\cos D \\
-\cos F & -\cos E & -\cos D & 1
\end{array} 
\right) .
$$
We here note that the following 
identity\footnote{Professors Derevnin and Mednykh thankfully mentioned to the author that 
a more general result was already known as Jacobi's theorem 
(see \cite[2.5.1.\ {\sc Theorem} (Jacobi)]{pr}).}
holds:
\begin{equation}
c_{3 \, 4}^2 - c_{3 \, 3} \, c_{4 \, 4} = \left( - \det G \right) \sin^2 A . \label{eqn_det}
\end{equation}

\vspace{\medskipamount}

We first consider the case where the endpoints of $A$ are both finite or ultraideal vertices.
In this case, since $c_{3 \, 3} \, c_{4 \, 4} > 0$, 
using (\ref{eqn_dist}),  (\ref{eqn_distance}) and (\ref{eqn_vertices}) we have
$$
- \cosh l = \left< \text{\boldmath $v$}_3, \text{\boldmath $v$}_4 \right>
= - \frac{c_{3 \, 4}}{\sqrt{c_{3 \, 3} \, c_{4 \, 4}}} ,
$$
where $l$ is the length of the edge joining $\text{\boldmath $v$}_3$ and $\text{\boldmath $v$}_4$.
Solve this equation by $\exp l$ and we have
$$
\exp l = \frac{c_{3 \, 4} + \sqrt{c_{3 \, 4}^2-c_{3 \, 3} \, c_{4 \, 4}}}{\sqrt{c_{3 \, 3} \, c_{4 \, 4}}} .
$$
Thus, using (\ref{eqn_det}), we have
\begin{equation}
\exp 2 l = \frac{2 c_{3 \, 4}^2 - c_{3 \, 3} \, c_{4 \, 4} + 2 c_{3 \, 4} \sqrt{-\det G} \sin A}
{c_{3 \, 3} \, c_{4 \, 4}} . \label{eqn_2l_1}
\end{equation}

\vspace{\smallskipamount}

We secondly consider the case where the endpoints of $A$ are a finite vertex and a ultraideal vertex.
In this case $c_{3 \, 3} \, c_{4 \, 4} <0$. So, using (\ref{eqn_distPoPl}), we have
$$
- \sinh l = \left< \text{\boldmath $v$}_3, \text{\boldmath $v$}_4 \right>
= - \frac{c_{3 \, 4}}{\sqrt{- c_{3 \, 3} \, c_{4 \, 4}}} .
$$
Solve this equation by $\exp l$ and, as (\ref{eqn_2l_1}), we have
\begin{eqnarray}
\lefteqn{\exp l = \frac{c_{3 \, 4} + \sqrt{c_{3 \, 4}^2 - c_{3 \, 3} \, c_{4 \, 4}}}
{\sqrt{- c_{3 \, 3} \, c_{4 \, 4}}}} \nonumber \\
&& \Longleftrightarrow \exp 2 l = \frac{2 c_{3 \, 4}^2 - c_{3 \, 3} \, c_{4 \, 4} + 2 c_{3 \, 4} \sqrt{-\det G} \sin A}
{- c_{3 \, 3} \, c_{4 \, 4}} . \label{eqn_2l_2}
\end{eqnarray}

\vspace{\medskipamount}

Next we calculate the partial derivative of ${\cal V} (T)$ with respect to $A$.
We first calculate $\partial U ( z_i, T ) / \partial A$;
\begin{eqnarray}
\frac{\partial U ( z_i, T )}{\partial A} 
&=& \frac{\partial U ( z_i, T )}{\partial a} \frac{\mathrm{d} a}{\mathrm{d} A} \nonumber\\
&=& \frac{\sqrt{-1}}{2} \log \frac{\varphi (z_i)}{\psi (z_i)} 
+ \frac{\partial U ( z_i, T )}{z_i} \frac{\partial z_i}{\partial a} \frac{\mathrm{d} a}{\mathrm{d} A},
\end{eqnarray}
where $\varphi (z) := \left( 1+ a b c z \right) \left( 1+ a e f z \right)$ and
$\psi (z) := \left( 1- a b d e z \right) \left( 1- a c d f z \right)$. 
It is easily checked that the second term of the right hand side is real, so we obtain
$$
\Im \frac{\partial U ( z_i, T )}{\partial A} = \frac{1}{2} \log \left| \frac{\varphi (z_i)}{\psi (z_i)} \right| ,
$$
and thus we have
\begin{eqnarray*}
\frac{\partial {\cal V} (T)}{\partial A}
&=& 
\frac{1}{4} \log \left| \frac{\varphi (z_1) \psi (z_2)}{\varphi (z_2) \psi (z_1)} \right| .
\end{eqnarray*}
By Theorem~\ref{thm_schlafli}, all we have to show is that 
$$
\frac{1}{4} \log \left| \frac{\varphi (z_1) \psi (z_2)}{\varphi (z_2) \psi (z_1)} \right| = -\frac{l}{2}
\Longleftrightarrow \left| \frac{\varphi (z_2) \psi (z_1)}{\varphi (z_1) \psi (z_2)} \right| = \exp{2 l} .
$$

\vspace{\medskipamount}

Substitute $z_1$ and $z_2$ for 
$\varphi (z_2) \psi (z_1) / \left( \varphi (z_1) \psi (z_2) \right)$ and we have
\begin{eqnarray*}
\frac{\varphi (z_2) \psi (z_1)}{\varphi (z_1) \psi (z_2)} 
&=& \frac{c_{3 \, 4} + \sqrt{-\det G} \sin{A}}{c_{3 \, 4} - \sqrt{-\det G} \sin{A}}\\
&=& \frac{2 c_{3 \, 4}^2 - c_{3 \, 3} c_{4 \, 4} + 2 c_{3 \, 4} 
	\sqrt{-\det G} \sin A}{c_{3 \, 3} c_{4 \, 4}} .
\end{eqnarray*}
Thus, as we saw in (\ref{eqn_2l_1}) and (\ref{eqn_2l_2}), 
the absolute value of the above is equal to $\exp{2 l}$.
Furthermore $\vol (T)$ and ${\cal V} (T)$ are 
analytic function whose variables are $A, B, \ldots , F$. \sloppy 
This means that ${\cal V} (T)$ is a volume formula for generalized tetrahedra up to some constant.
Now we check that
${\cal V} (T)$ is actually the formula.

\vspace{\bigskipamount}

Consider the regular tetrahedron in ${\Bbb P}^3$ circumscribing $\overline{{\bf B}^3}$.
Then the polyhedron obtained by the truncation at all vertices 
is the regular ideal octahedron with all dihedral angles being $\pi/2$,
which we denote by $T$.
Its volume
is calculated as follows:
since $T$ is symmetric, it can be carved like the fringe of an orange and becomes
four isometric ideal tetrahedra of the triple of dihedral angles are $\pi/2$, $\pi/4$ and $\pi/4$.
Using Milnor's volume formula for ideal tetrahedra (see \cite{mi}) we have
\begin{eqnarray*}
\vol (T) &=& 4 \left( \Lambda (\frac{\pi}{2}) + \Lambda (\frac{\pi}{4}) + \Lambda (\frac{\pi}{4}) \right) \\
&=& 8 \Lambda (\frac{\pi}{4}) \left( \approx 3.66386 \right) .
\end{eqnarray*}
Here we note that the Lobachevsky function $\Lambda (x)$
is a periodic odd function with period $\pi$ and $\Lambda (\pi/2) = 0$.
Furthermore there is the following relationship 
with the dilogarithm function (see, for example, \cite[\S 1.1.4]{ki})
\begin{equation}
\Im \li (\exp \sqrt{-1} x) = 2 \Lambda (\frac{x}{2}) \quad \text{for any $x \in {\Bbb R}$.}
\label{eqn_lob}
\end{equation}

On the other hand, $T$ is regarded as the truncated tetrahedron of all dihedral angles being $0$.
So we have
$$
(z_1 , z_2) = ( -\sqrt{-1}, \sqrt{-1} ) .
$$
Moreover in this case 
$$
U ( z, T ) = 2 \left( \li (z) - \li (- z) \right) .
$$
So, using (\ref{eqn_lob}), we have
\begin{eqnarray*}
{\cal V} (T) &=& \Im \left(\li (z_2) - \li (-z_2) - \li (z_1) + \li (z_1) \right) \\
&=& 2 \Im \left(\li (\sqrt{-1}) - \li (-\sqrt{-1}) \right)\\
&=& 4 \left(\Lambda (\frac{\pi}{4}) - \Lambda (\frac{-\pi}{4})) \right)\\
&=& 8 \Lambda (\frac{\pi}{4})\\
&=& \vol (T) ,
\end{eqnarray*}
which show that ${\cal V} (T)$ is actually the volume of generalized tetrahedra $T$.
Thus we have proved Theorem~\ref{thm_volumeformula}. \hfill $\square$

\chapappendix{Murakami-Yano's volume formula}
\label{ap_MY}

Murakami-Yano's volume formula comes from the quantum $6j$-symbol, 
which is a source of invariants for three-dimensional manifolds. 
To obtain the formula they applied R.~ Kashaev's method of computation the quantum $6j$-symbol.
He discovered in \cite{ka} a relationship between certain quantum invariants 
and the volumes of several hyperbolic knots in the three-dimensional sphere.
A technical key tool in his paper is a {\em saddle point approximation}, 
which is an asymptotic approximation of some integral by saddle point of an integrand.
The values $z_1$ and $z_2$ in this paper corresponds to saddle points of $U ( z, T )$,
namely, roughly speaking, $z_1$ and $z_2$ 
arise from the solution of $\mathrm{d} U (z,T)/ \mathrm{d} z = 0$.

They also proved in \cite{my} that, by taking an appropriate branch of $U ( z, T )$, 
we have $\Re {\cal V} (T) = 0$ if $T$ is a hyperbolic tetrahedron 
and ${\cal V} (T) = \vol (T)$ if $T$ is a spherical tetrahedron.

\chapappendix{Graphs on the volume and the edge length of a regular tetrahedron}
\label{ap_graph}

The following graphs show the correspondences
between the dihedral angle and the volume or the edge length of regular tetrahedra.
The dihedral angle varies from $0$ to 
$\arccos (1/3) \approx 1.230959$.
The angle $0$ corresponds to the regular ideal octahedron as we saw in Theorem~\ref{thm_max},
and the angle 
$\arccos (1/3)$
corresponds to the Euclidean regular tetrahedron,
corresponding to infinitesimally small  hyperbolic regular tetrahedra.
The dihedral angle of the regular ideal  tetrahedron is $\pi/3$.
\begin{figure}[ht]
        \begin{center}
        \includegraphics[clip]{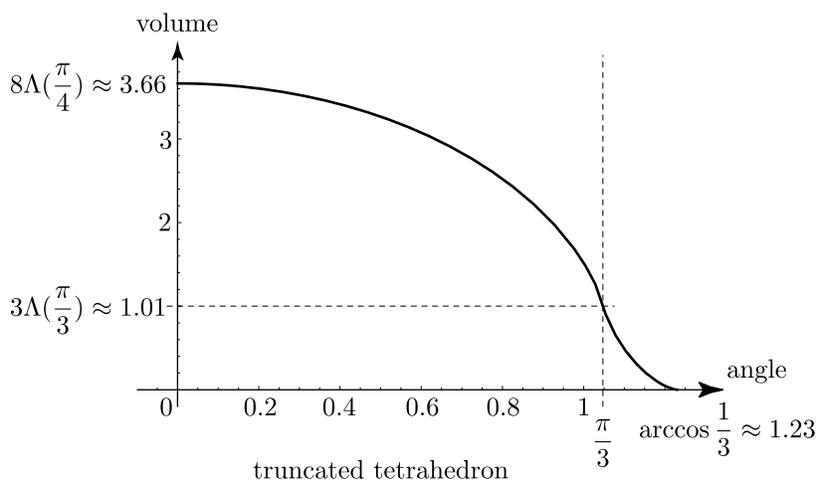}  
        \end{center}
        \caption{Volume of a regular tetrahedron}
        \label{fig_volume}
\end{figure}
\begin{figure}[ht]
        \begin{center}
        \includegraphics[clip]{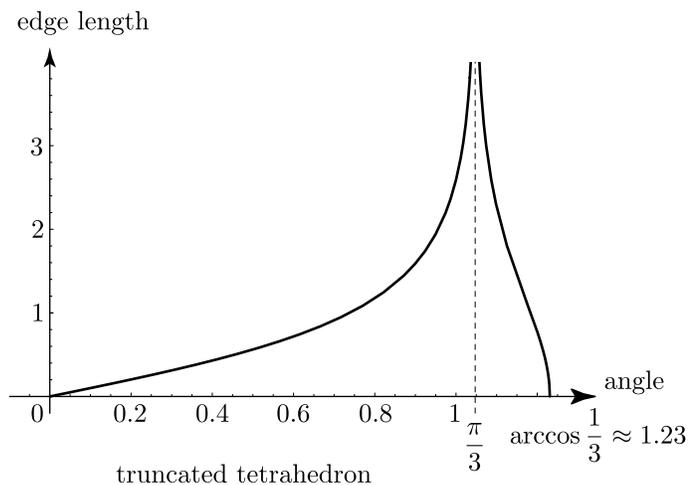}  
        \end{center}
        \caption{Edge length of a regular tetrahedron}
        \label{fig_edge}
\end{figure}

\chapappendix{Further generalization and open problems; when the truncation is not ``mild"}
\label{ap_Lambert}

In this paper it is assumed that the truncation of a tetrahedron is ``mild," namely 
each faces obtained by the truncations do not intersect each other in ${\Bbb H}^3$
(this is equivalent to the assumption of the interior of all edges intersecting $\overline{{\bf B}^3}$;
see Definition~\ref{df_trunc}).
What happens if we do not assume this condition, i.e., if the truncation is ``intense"?

A three-dimensional {\em orthoscheme\/} is a tetrahedron 
with dihedral angles $B=E=F=\pi/2$ (see \cite{ke} for details; 
another basic reference is \cite{bh}\footnote{Professor Moln\'ar 
thankfully told the author about this book.}).
Then the edge spanned by 
$\text{\boldmath $v$}_1$ and $\text{\boldmath $v$}_2$
is perpendicular to the face spanned by 
$\text{\boldmath $v$}_2 , \text{\boldmath $v$}_3 , \text{\boldmath $v$}_4$, 
and also the face spanned by 
$\text{\boldmath $v$}_1 , \text{\boldmath $v$}_2 , \text{\boldmath $v$}_3$
is perpendicular to the edge spanned by 
$\text{\boldmath $v$}_3$ and $\text{\boldmath $v$}_4$,
which is the meaning of 
``An {\em orthoscheme\/} is a simplex of a generalization (with respect to dimensions) of a right triangle"
mentioned in Introduction.
Two vertices $\text{\boldmath $v$}_1$ and $\text{\boldmath $v$}_4$ 
are called {\em principal vertices\/}.

For an orthoscheme, 
if two faces obtained by truncations at principal vertices intersect each other in the hyperbolic space,
then the polyhedron becomes a cube with all faces Lambert quadrilaterals 
(i.e., quadrilaterals with one acute and three right angles).
Such a cube is called {\em Lambert cube\/}.
In \cite{ke} she also constructed a volume formula for Lambert cubes.
This formula is closely related to the one for orthoschemes,
so the following problem naturally arises:
%
%
%
%
%
%
%
%
\begin{pbm}
Construct a volume formula for intensely truncated tetrahedra from our volume formula.
\end{pbm}

One of the first step of solving this problem is to see 
the relationship between the formula for mildly truncated orthoschemes in \cite{ke}
and our formula.
Her formula is presented via Lobachevsky functions, 
and we can easily rewrite our formula with respect to it (see \cite[$(1.19)$]{my}).
Theoretically ours contains hers, but there is no direct transformation.
So the following problem also arises:
%
%
%
%
%
%
%
%
\begin{pbm}
For orthoschemes, find a transformation between our volume formula and Kellerhals's formula.
\end{pbm}
%


\end{document}